\documentclass[11pt]{amsart}

\usepackage{mathrsfs}
\usepackage{amssymb}

\usepackage{titletoc}
\input xy
  \xyoption{all}
\usepackage{amscd}
\usepackage{amsmath, amssymb}
\usepackage{amsfonts}
\usepackage[colorlinks,linkcolor=blue,citecolor=blue, pdfstartview=FitH]{hyperref}

\numberwithin{equation}{section}

\theoremstyle{plain}
\newtheorem{thm}{Theorem}[section]
\newtheorem{theorem}[thm]{Theorem}
\newtheorem{lemma}[thm]{Lemma}
\newtheorem{corollary}[thm]{Corollary}
\newtheorem{proposition}[thm]{Proposition}

\theoremstyle{definition}

\newtheorem{remark}[thm]{Remark}

\newtheorem{definition}[thm]{Definition}

\newtheorem{example}[thm]{Example}

\newtheorem{defn-thm}[thm]{Definition-Theorem}

\newcommand{\sF}{{\mathcal F}}

\newcommand{\sI}{{\mathcal I}}

\newcommand{\sO}{{\mathcal O}}

\newcommand{\sX}{{\mathcal X}}

\newcommand{\C}{{\mathbb C}}

\newcommand{\K}{{\mathbb K}}

\renewcommand{\P}{{\mathbb P}}

\newcommand{\R}{{\mathbb R}}

\newcommand{\Z}{{\mathbb Z}}

\newcommand{\ch}{{ ch}}
\newcommand{\qtq}[1]{\quad\mbox{#1}\quad}
\newcommand{\bp}{\bar{\partial}}

\newcommand{\ds}{\oplus}

\newcommand{\ts}{\otimes}

\newcommand{\btheorem}{\begin{theorem}}
\newcommand{\etheorem}{\end{theorem}}
\newcommand{\bproposition}{\begin{proposition}}
\newcommand{\eproposition}{\end{proposition}}
\newcommand{\bdefinition}{\begin{definition}}
\newcommand{\edefinition}{\end{definition}}
\newcommand{\bcorollary}{\begin{corollary}}
\newcommand{\ecorollary}{\end{corollary}}
\newcommand{\bproof}{\begin{proof}}
\newcommand{\eproof}{\end{proof}}
\newcommand{\bremark}{\begin{remark}}
\newcommand{\eremark}{\end{remark}}
\newcommand{\eexample}{\end{example}}
\newcommand{\bexample}{\begin{example}}
\newcommand{\la}{\langle}
\newcommand{\elemma}{\end{lemma}}
\newcommand{\blemma}{\begin{lemma}}
\newcommand{\ra}{\rangle}
\newcommand{\sq}{\sqrt{-1}}

\newcommand{\p}{\partial}

\renewcommand{\bar}{\overline}
\newcommand{\eps}{\varepsilon}

\renewcommand{\phi}{\varphi}

\newcommand{\ee}{\end{eqnarray*}}
\newcommand{\be}{\begin{eqnarray*}}

\newcommand{\beq}{\begin{equation}}
\newcommand{\eeq}{\end{equation}}

\newcommand{\bd}{\begin{enumerate}}
\newcommand{\ed}{\end{enumerate}}

\renewcommand{\tilde}{\widetilde}

\renewcommand{\>}{\rightarrow}

\usepackage{fancyhdr}
\renewcommand{\leftmark}{{\scriptsize Global generation and very ampleness for adjoint linear series}}

\begin{document}
\title{Global generation and very ampleness for adjoint linear series }
\makeatletter
\let\uppercasenonmath\@gobble% disables title uppercase
\let\MakeUppercase\relax% disables author uppercase
\let\scshape\relax% disables section smallcaps
\makeatother
\author{Xiaoyu Su and Xiaokui Yang}
%\date{ }

\address{{Address of Xiaoyu Su:
        Academy of Mathematics and Systems Science,
        Chinese Academy of Sciences, Beijing, 100190, China.}}
\email{\href{mailto:suxiaoyu14@ucas.ac.cn}{{suxiaoyu14@mails.ucas.ac.cn}}}
\address{{Address of Xiaokui Yang: Morningside Center of Mathematics, Institute of
        Mathematics; Hua Loo-Keng Center for Mathematical Sciences,
        Academy of Mathematics and Systems Science,
        Chinese Academy of Sciences, Beijing, 100190, China.}}
\email{\href{mailto:xkyang@amss.ac.cn}{{xkyang@amss.ac.cn}}}
\maketitle

\begin{abstract} Let $X$ be a smooth projective variety over
an algebraically closed field $\K$ with arbitrary characteristic.
Suppose $L$ is an ample and globally generated line bundle. By
Castelnuovo--Mumford regularity, we show that
 $K_X \ts L^{\ts \dim X} \otimes A$ is globally generated and  $K_X \ts L^{\ts (\dim X+1)} \otimes A$ is very ample,
  provided the line bundle $A$ is  nef but not  numerically trivial.
On complex projective varieties, by investigating
Kawamata-Viehweg-Nadel type vanishing theorems for vector bundles,
we also obtain the global generation for adjoint vector bundles. In
particular, for a holomorphic submersion $f:X\>Y$ with   $L$ ample
and globally generated, and $A$ nef but not numerically trivial, we
 prove the global generation of
$ f_*(K_{X/Y})^{\ts s}\ts K_Y \ts L^{\ts \dim Y} \otimes A$ for any
positive integer $s$.

\end{abstract}

\small{\setcounter{tocdepth}{1} \tableofcontents}

\section{Introduction}

    In \cite{Fuj88}, Fujita proposed the
   following conjecture\\
\noindent\textbf{Conjecture.}  Let $X$ be a smooth projective
variety
   and $L$ be an ample line bundle. Then
\bd \item $K_X \ts L^{\ts (\dim X+1)} $ is globally generated;

\item $K_X \ts L^{\ts (\dim X+2)}$ is very ample.

\ed

\noindent Fujita's conjecture is a deceptively simple open question
in classical algebraic geometry. Up to dimension 4, the global
generation conjecture has been proved (\cite{Rei88, EL93, Kaw97}).
Recently,  Fei Ye and Zhixian Zhu proved  it in dimension $5$
(\cite{YZ2}).  Also, there are many other ``Fujita Conjecture type"
theorems have been proved, and we refer the reader to
 \cite{Dem93,  Kol93, AS95, Siu96, Sm97, Hel97,  Hel99, Sm00, Hei02, Ara04, Kee08, PP08, BBMT14, PS14, Sch14, YZ1, Fuj, FS}) and the references therein.

 In this paper, we prove Fujita Conjecture type theorems by using
  analytic methods in complex geometry when the twisted line (resp. vector) bundle  is nef (resp. Nakano semi-positive).

\subsection{Fujita Conjecture type theorems on projective varieties over arbitrary fields}
Let $\K$ be an algebraically closed field with arbitrary
characteristic. By using characteristic $p$ methods, Keeler proved
in \cite[Theorem~1.1]{Kee08} the following interesting result.

\btheorem[Keeler]\label{keeler} Let $X$ be a smooth projective
variety over $\K$ with dimension $n$. Suppose $L$ is an ample and
globally generated line bundle, and $A$ is an ample line bundle.
Then

\bd \item $K_X \ts L^{\ts n} \otimes A$ is globally generated;

\item $K_X \ts L^{\ts (n+1)} \otimes A$ is very ample.

\ed

\etheorem

\noindent   In the case $\K=\C$, Angehrn and Siu  proved the very
ampleness part of Theorem \ref{keeler}  by analytic methods
(\cite[Lemma~11.1]{AS95}).  In \cite[Corollary~4.8]{Sch14}, Schwede
generalized the global generation part of Theorem \ref{keeler} to
the case when $A$ is nef and big, and $\ch(\K)=p>0$. The first
result of our paper deals with a more general case when $A$ is nef
but not numerically trivial.

\btheorem\label{main1}  Let $X$ be a smooth projective variety over
$\K$ with dimension $n$. Suppose $L$ is an ample and globally
generated line bundle.  If $A$ is a  nef but not numerically trivial
line bundle, then

\bd \item $K_X \ts L^{\ts n} \otimes A$ is globally generated;

\item $K_X \ts L^{\ts (n+1)} \otimes A$ is very ample.

\ed \etheorem

\noindent As far as the authors know, this is the greatest
generality in which Fujita  Conjecture type theorem has been proved
in any characteristic. Theorem \ref{main1} is also optimal in the
sense that one can not drop the non-triviality condition on $A$,
which can be seen from the  example $(X,L)=(\P^n,\sO(1))$. One may
want to know the limit case when $A$ is indeed a trivial line
bundle. For the global generation part, Smith proved in
\cite[Theorem~2]{Sm00} that, this example is the only exceptional
case, i.e.,  if $(X,L)\neq (\P^n,\sO(1))$, then $K_X\ts L^{\ts n}$
is globally generated. Her proof relies on the ``tight closure"
methods in the frame of commutative algebra.

\subsection{Fujita Conjecture type theorems on complex projective varieties}

In this subsection, we  focus on the cases in complex geometry, i.e.
$X$ is defined over the complex number field $\C$. At first, we
obtain a  general version of the global generation part in Theorem
\ref{main1}:

\btheorem\label{main11}Let $X$ be a compact K\"ahler manifold of
dimension $n$ and $L$ be an ample and globally generated line
bundle. Suppose $(A,e^{-2\phi})$ be a pseudo-effective line bundle
and $\sI(\phi)$ is the multiplier ideal sheaf. If the numerical
dimension of  $(A,\phi)$ is not zero, i.e. $\emph{nd}(A,\phi)\neq
0$, then
 \beq K_X \ts L^{\ts n} \otimes A\ts \sI(\phi)\eeq is globally
generated. \etheorem

\noindent A key ingredient in the proof of Theorem \ref{main11}
relies on  vanishing theorems for pseudo-effective line bundles
(\cite[Corollary~1.7]{GZ}, \cite[Theorem~0.15]{Dem14},
\cite[Corollary~3.2]{GZ15} or a weaker version
\cite[Theorem~1.3]{Cao14}).

By using analytic methods, we  also investigate the globally
generated property for adjoint vector bundles.

\btheorem\label{main2} Let $X$ be a compact K\"ahler manifold of
dimension $n$ and $L$ be an ample and globally generated line
bundle. Let $(E,h)$ be a Hermitian  holomorphic vector bundle with
Nakano semi-positive curvature. Suppose $A$ is a nef but not
numerically trivial line bundle, then the adjoint vector bundle \beq
K_X \ts L^{\ts n} \otimes (E\ts A)\eeq is globally generated.
\etheorem

\noindent Theorem \ref{main2} is derived from the following
vanishing theorem for vector bundles, building on  ideas in
\cite{GZ}, \cite{Dem} and \cite{Cao14}.

\btheorem\label{mainkey} Let $X$ be a  complex projective variety
with $\dim X=n$. If $(A,e^{-2\phi})$ is a pseudo-effective line
bundle and $(E,h)$ a Nakano semi-positive vector bundle, then \beq
H^q(X,K_X\ts E\ts A\ts \sI(\phi))=0 \qtq{for} q>n-
\emph{nd}(A,\phi).\eeq \etheorem

\noindent\textbf{Remark.}  According to \cite{GZ2}, it is not hard
to see that Theorem \ref{mainkey} is also true on compact K\"ahler
manifolds. Hence, there is a version of Theorem \ref{main2} for
pseudo-effective line bundle $(A,e^{-2\phi})$ (e.g. Theorem
\ref{main211}). For simplicity, we only formulate applications for
\emph{ nef line bundles} (see Theorem \ref{main311} for general
cases).

 As an application of Theorem \ref{main2}, we obtain  the
following result in pure algebraic language.

\btheorem\label{main3} Let $f:X\rightarrow Y$ be a holomorphic
submersion between two  complex projective varieties and $\dim Y=n$.
Suppose $L\>Y$ is an ample and globally generated line bundle, and
$A\>Y$ is a nef but not numerically trivial line bundle. Then \beq
f_*(K_{X/Y})^{\ts s}\ts K_Y \ts L^{\ts n} \otimes A\eeq is globally
generated for any $s\geq 1$. \etheorem

As a special case of Theorem \ref{main3}, we obtain the following
well-known result of Koll{\'a}r (\cite[Theorem~3.5,
Theorem~3.6]{Kol86}):

 \bcorollary[Koll{\'a}r] Let $f:X\>Y$ be a holomorphic
submersion between two smooth projective varieties and $\dim_\C
Y=n$. Suppose $L\>Y$ is an ample and globally generated line bundle,
then
 $$f_*(K_{X/Y})^{\ts s}\ts K_Y \ts L^{\ts
(n+1)}$$ is globally generated for $s\geq 1$. \ecorollary

 As another application of Theorem \ref{main2}, we also get global generation of pluricanonical adjoint
bundles of canonically polarized families.

\btheorem\label{main4} Let $f:X\>S$ be a holomorphic family of
canonically polarized compact K\"ahler manifolds  effectively
parameterized by a smooth projective variety $S$ with dimension $n$.
Suppose $L\>S$ is an ample and globally generated line bundle, and
$A\>S$ is a nef line bundle. Then \beq f_*(K^{\ts s}_X) \ts L^{\ts
n} \otimes A\eeq is globally generated for $s>1$. \etheorem

\vskip 2\baselineskip

\section{Fujita Conjecture type theorems on projective varieties over arbitrary fields}
In this section, we investigate Fujita Conjecture type theorems on
projective varieties over algebraically closed fields and prove
Theorem \ref{main1}. Let $\K$ be an algebraically closed field with
arbitrary characteristic, and $X$ be a smooth projective variety
over $\K$. Firstly we introduce the theory of Castelnuovo--Mumford
regularity.

\subsection{Castelnuovo--Mumford regularity}

 Suppose $L$ is an ample and globally generated
line bundle over $X$.

\bdefinition\label{cm} A coherent sheaf $\sF$ on $X$ is $m$-regular
with respect to $L$ if \beq H^q(X,\sF \ts L^{\ts (m-q)})=0 \qtq{for
} q>0.\eeq \edefinition

\noindent The following results are well-known (e.g.
\cite[Section~1.8]{Laz04}, or \cite[Section~5.2]{FGIKNV05}), and for
the sake of completeness we include a proof here.
 \blemma[Mumford]\label{cmr} Let $\sF$ be a $0$-regular coherent sheaf on
$X$ with respect to $L$, then $\sF$ is generated by its global
sections.  \elemma

\bproof Suppose $\dim X=n$. We shall use standard hyperplane
induction method to prove it. For simplicity, we write $L=\sO_X(1)$.
Since the coherent sheaf $\sF$ has finitely many associated points,
we can choose $s\in H^0(X,L)$ such that the corresponding divisor
$B$ does not contain any associated point of $\sF$. Hence, for any
$i\geq 0$, we have the exact sequence \beq 0\>\sF(-i-1)\stackrel{\ts
s}{\>}\sF(-i)\>\sF_B(-i)\>0.\label{a}\eeq By using the associated
long exact sequence \beq
\cdots\>H^{i}(X,\sF(-i))\>H^{i}(X,\sF_B(-i))\>H^{i+1}(X,\sF(-i-1))\>\cdots\eeq
and $0$-regularity of $\sF$, we see $ H^i(X,\sF_B(-i))=0$ for $i>0$,
i.e. $\sF_B$ is $0$-regular with respect to $L$. Similarly, for any
$q\geq 0$, we have \beq \cdots
\>H^q(X,\sF(-q))\>H^{q}(X,\sF(1)(-q))\>H^{q}(X,\sF_B(1)(-q))\>\cdots.
\label{re}\eeq We show $\sF(1)$ is $0$-regular if $\sF$ is
$0$-regular. Indeed, we have $\sF_B$ is $0$-regular and by
hyperplane induction hypothesis, $\sF_B(1)$ is $0$-regular. By
(\ref{re}),   $\sF(1)$ is $0$-regular. Hence, we know $\sF(k)$ is
$0$-regular for all $k\geq 0$. In the commutative diagram
$$\xymatrix{
  & H^0(X,\sF(k))\ts H^0(X,\sO(1)) \ar[d]^{\iota_X} \ar[r]^{r_k\ts 1}
                      & H^0(X,\sF_B(k))\ts H^0(X,\sO_X(1))  \ar[d]_{\iota_B}    \\
  H^0(X,\sF(k)) \ar[r]^{\ts s}& H^0(X,\sF(k+1)) \ar[r]^{r_{k+1}}     & H^0(X,\sF_B(k+1)),               }
$$
$r_{k}\ts 1$ and $r_{k+1}$ are surjective according to long exact
sequence associated to (\ref{a}). By diagram chasing, it is obvious
that $\iota_X$ is surjective if and only if $\iota_B$ is surjective.
Note that in the above commutative diagram we can replace $B$ by
intersections of suitable divisors in $|L|$. Hence, we can show
$\iota_X$ is surjective by  induction on $\dim B$. When $\dim B=0$,
it is easy to see $\iota_B$ is surjective. Hence, by induction, we
know $\iota_X$ is surjective. For any $x\in X$ and large $N$, in the
commutative diagram
$$\xymatrix{
   H^0(X,\sF)\ts H^0(X,\sO(1))^{\ts N} \ar[d]^{1\ts ev_1
   } \ar[r]^{\ \ \ \ \ \ \ \ \ \ \ \iota}
                      & H^0(X,\sF(N)) \ar[d]_{ev_2}    \\
  H^0(X,\sF)\ts (\sO(1)|_x)^{\ts N}  \ar[r]^{\ \ \ \ \ \ \ \ f}     & (\sF(N))|_x,               }
$$
$\iota$ is surjective since $\iota_X$ is surjective. On the other
hand, $1\ts ev_1$ and $ev_2$ are both surjective and so $f$ is
surjective, and we deduce $H^0(X,\sF)\>\sF|_x$ is surjective.
 \eproof

\subsection{The proof of Theorem \ref{main1}} Before giving the proof of Theorem \ref{main1}, we present a more general result.

\btheorem\label{variant} Let $X$ be a smooth projective variety over
 $\K$  with dimension $n$. Suppose $L$ is an ample
and globally generated line bundle. If $A$ is a line bundle
 such that $L\ts A$ is ample and the Kodaira-Iitaka dimension
 $\kappa(A^*)=-\infty$, then

\bd \item $K_X \ts L^{\ts n} \otimes A$ is globally generated;

\item $K_X \ts L^{\ts (n+1)} \otimes A$ is very ample.

\ed

\etheorem

\bproof Suppose  $ch(\K)=p>0$.   We shall show that $F_*^k(K_X)\ts
L^{\ts n} \ts A$ is $0$-regular for all $k\geq k_0$ where $F:X\>X$
is the absolute Frobenius morphism. When $0<q<n$, i.e. $n-q\geq 1$,
$L^{\ts(n-q)}\ts A$ is ample since both $L$ and $L\ts A$ are ample.
Hence, by the Serre vanishing theorem, for each $q$ with $0<q<n$,
there exists a positive constant $k_q=k(q)>0$ such that \beq
H^q(X,K_X \ts L^{\ts p^k(n-q)}\ts A^{\ts p^k} )=0,\eeq for $k\geq
k_q$. By the projection formula, one has \beq F_*^k(K_X \ts L^{\ts
p^k(n-q)}\ts A^{\ts p^k})=F_*^k(K_X)\ts L^{\ts (n-q)}\ts A.\eeq
Since the Frobenius morphism $F$ is a finite morphism, we get \be
H^q(X,F_*^k(K_X)\ts L^{\ts (n-q)} \ts A)&\cong &H^q(X,F_*^k(K_X \ts
L^{\ts p^k(n-q)}\ts A^{\ts p^k}))\\
&\cong &H^q(X,K_X \ts L^{\ts
p^k(n-q)}\ts A^{\ts p^k} )\\
&=&0.\ee

\noindent  When $q=n$, we want to show \beq  H^n(X, F_*^k(K_X)\ts
A)=0.\label{n}\eeq
    By the projection formula again, we  have
    $  F_*^k(K_X) \ts A = F_*^k\left(K_X \ts A^{\ts p^k}\right)
    $
    and  \be H^n(X, F_*^k(K_X)\ts A)&\cong &H^n\left(X,
F_*^k(K_X \ts
    A^{\ts p^k})\right)\\&\cong &H^n(X, K_X \ts
    A^{\ts p^k})\\&\cong & H^0(X,(A^{\ts p^k})^*)^*. \ee  Since the Kodaira-Iitaka dimension
    $\kappa(A^*)=-\infty$, we have
    $  H^0(X,(A^*)^{\ts \ell})=0$  for any $\ell>0$. Hence we get
    (\ref{n}).  By Definition \ref{cm},  $F_*^k(K_X)\ts A\ts L^{\ts n}$ is
$0$-regular for all $k\geq k_0$ where $k_0=\max\{k_1,\cdots,
k_{n-1}\}$. According to Lemma \ref{cmr},  $F_*^k(K_X)\ts L^{\ts n}
\ts A$ is globally
   generated.  Thanks to \cite[Lemma~3.3]{Kee08}, $K_X \ts L^{\ts n} \ts A$ is a quotient of
$$F_*^k(K_X)\ts A\ts L^{\ts n}$$ for all $k\geq k_0$ and so $K_X \ts
L^{\ts n} \ts A$ is globally generated.

  It is well-known that  $K_X\ts L^{\ts (n+1)}\ts A$ is very ample
  if and only if  for every $x\in X$, $\mathfrak{m}_x\ts K_X\ts L^{\ts(n+1)}\ts
  A$ is globally generated. Since $ F_*^k(K_X)\ts L^{\ts n}\ts A$ is
  $0$-regular with respect to $L$, it is proved in
  \cite{Kee08} that for every $x\in X$, $\mathfrak{m}_x
\ts F_*^k(K_X)\ts L^{\ts (n+1)}\ts A$ is also $0$-regular and so it
is globally generated. By \cite[Lemma~3.3]{Kee08},  $\mathfrak{m}_x
\ts K_X \ts L^{\ts (n+1)}\ts A$ is a quotient  of $\mathfrak{m}_x
\ts F_*^k(K_X)\ts L^{\ts (n+1)}\ts A$. Hence $\mathfrak{m}_x\ts
K_X\ts L^{\ts(n+1)}\ts
  A$ is  globally generated.

 When $ch(\K)=0$, $L^{\ts
k}\ts A$ is ample for every $k>0$. Hence by Kodaira vanishing
theorem, we have
 \beq H^q(X,K_X\ts L^{\ts (n-q)}\ts A)=0 \eeq
 for $0<q<n$. Since $\kappa(A^*)=-\infty$, we also have  \beq H^n(X,K_X\ts A)\cong
 H^0(X,A^*)=0.\eeq Hence,  by Lemma \ref{cmr} again, we know  $K_X\ts L^{\ts n}
\ts A$ is $0$-regular with respect to $L$ and hence it is globally
generated. The very ampleness of  $K_X\ts L^{\ts (n+1)} \ts A$ can
be proved similarly.
 \eproof

\blemma\label{Lemma01} Let $A$ be a nef but not numerically trivial
line bundle over a smooth projective variety $X$ over $\mathbb{K}$.
Then the Kodaira-Iitaka dimension $\kappa(A^*)=-\infty$, i.e.
$H^0(X, (A^*)^{\ts \ell})=0$ for all $\ell>0$. \elemma

\bproof  We also use $A$ to denote the divisor class of $A$.
 Since $A$ is nef but not numerically trivial, it is well-known that (e.g. \cite[Section~3.8]{Deb13}) that  there exists an ample divisor $H$,
 such that $A \cdot H^{n-1}> 0$. We show $\kappa(A^*)=-\infty$.   Suppose $H^0(X,\sO(-\ell A)) \neq 0$ for some
$\ell>0$. Then $-\ell A$ is effective. Let $-\ell A=\sum \nu_i D_i$
where $D_i$ are irreducible divisors and $\nu_i\geq 0$. Since $-\ell
A$ is not numerically trivial, there is at least one $\nu_i>0$.   By
Nakai-Moishezon criterion for ampleness, for ample divisor $H$ in
$X$, we have $ H^{n-1}\cdot (-\ell A)>0$, i.e. $H^{n-1}\cdot A<0$
which contradicts $A\cdot H^{n-1}>0$. \eproof

\noindent \emph{The proof of Theorem \ref{main1}}. It follows from
Lemma \ref{Lemma01} and Theorem \ref{variant}.

\bremark Theorem \ref{variant} is more general than Theorem
\ref{main1}. In Theorem \ref{variant}, $A$ can be certain
numerically trivial line bundles.  It is easy to see that, if $A$ is
a non-torsion point in $Pic^0(X)$, Theorem \ref{variant} also works.
For instance, let $E=\C/(\Z\ds\sq \Z)$ be an elliptic curve. Suppose
$A=\sO(P-Q)$ where $P$ is a rational point and $Q$ is an irrational
point on $E$. Then $A$ is numerically trivial and
$\kappa(A^*)=-\infty.$  Indeed, $(A^*)^{\ts \ell}$ has no nonzero
section for any $\ell>0$. Otherwise, the divisor $\ell(P-Q)$ is
linearly equivalent to the zero divisor, which is absurd. \eremark

\vskip 2\baselineskip

\section{Vanishing theorems for vector bundles on compact K\"ahler manifolds}
\noindent In this section, we investigate various vanishing theorems
for vector bundles, which are the key ingredients in the proof of
Fujita Conjecture type theorems. In particular, we give the proof of
Theorem \ref{mainkey}.

Let $E$ be a holomorphic vector bundle over a compact complex
manifold $X$ and $h$ be a smooth Hermitian metric on $E$. There
exists a unique connection $\nabla$ which is compatible with the
Hermitian
 metric $h$ and the complex structure on $E$. It is called the Chern connection of $(E,h)$. Let $\{z^i\}_{i=1}^n$ be  local holomorphic coordinates
  on $X$ and  $\{e_\alpha\}_{\alpha=1}^r$ be a local frame
 of $E$. The curvature tensor $\Theta^E\in \Gamma(X,\Lambda^2T^*X\ts E^*\ts E)$ has components \beq R_{i\bar j\alpha\bar\beta}= -\frac{\p^2
h_{\alpha\bar \beta}}{\p z^i\p\bar z^j}+h^{\gamma\bar
\delta}\frac{\p h_{\alpha \bar \delta}}{\p z^i}\frac{\p
h_{\gamma\bar\beta}}{\p \bar z^j}.\eeq Here and henceforth we
sometimes adopt the Einstein convention for summation.

\bdefinition
A Hermitian holomorphic vector bundle $(E,h)$ is called Nakano
positive (resp. Nakano semi-positive) if
$$ R_{i\bar j
\alpha\bar\beta} u^{i\alpha}\bar u^{j\beta}>0\qtq{(resp. $\geq 0$)}
$$
for nonzero vector $u=(u^{i\alpha})\in \C^{nr}$. \edefinition

\noindent Let's describe some elementary properties on positive
vector bundles.

\blemma[Nakano vanishing theorem]\label{nakano} Let $X$ be a compact
K\"ahler manifold. Suppose $(E,h)\>X$ is a Hermitian holomorphic
vector bundle with Nakano positive curvature, then \beq H^q(X,K_X\ts
E)=0, \ \ q\geq 1.\eeq

\elemma

\blemma\label{positive} Let $(X,\omega_g)$ be a compact K\"ahler manifold. \bd
\item Let $(E,h)$ be a Nakano positive vector bundle and
$A$ a nef line bundle. Then $E\ts A$ admits a Hermitian metric with
Nakano positive curvature.

\item  Let $(E,h)$ be a Nakano semi-positive vector bundle and
$A$ be an ample line bundle. Then $E\ts A$ admits a Hermitian metric
with Nakano positive curvature.

\item Let $(E,h^E)$ and $(\tilde E,  h^{\tilde E})$ be two Nakano
semi-positive vector bundles, then $(E\ts \tilde E, h\ts\tilde h)$
is also Nakano semi-positive.

 \ed \elemma

\bproof (1).  For the fixed K\"ahler metric $\omega_g$ on $X$, there
exists a constant $\eps>0$ such that
$$\sq \Theta^E(u(x),u(x))\geq 2\eps |u(x)|_{g\ts h}^2$$
for any $u\in\Gamma(X,T^{1,0}X\ts E)$. On the other hand, by
analytic definition of nefness (e.g. \cite{Dem}), there exists a
smooth metric $h_0$ on the nef line bundle $A$ such that \beq \sq
\Theta^A \geq -\eps \omega_g. \eeq The curvature of $h\ts h_0$ on
$E\ts A$ is $ \Theta^{E\ts A}=\Theta^E\cdot id_A+ id_E\cdot
\Theta^A.$ Hence, for any $u\in \Gamma(X,T^{1,0}X\ts E)$ and $v\in
\Gamma(X,A)$ \beq \sq\Theta^{E\ts A}(u\ts v,u\ts v)\geq
\left(\sq\Theta^E(u,u)-\eps |u|_{g\ts h}^2\right) |v|^2_{h_0}\geq
\eps |u|_{g\ts h}^2|v|^2_{h_0}. \eeq Therefore, $E\ts A$ is Nakano
positive. The proof of $(2)$ is similar to that of $(1)$.

(3). By using  curvature formula of $h\ts \tilde h$ on $E\ts \tilde
E$, $$ \Theta^{E\ts A}=\Theta^E\cdot id_{\tilde E}+ id_E\cdot
\Theta^{\tilde E},$$  for any (local) vectors $u\in
\Gamma(X,T^{1,0}X\ts E \ts \tilde E)$ with the form $u=u^{i\alpha
A}dz^i\ts e^\alpha\ts e^A$ in the local holomorphic frames $\{z^i,
e^\alpha, e^A\}$ of $\{X, E,\tilde E\}$, we obtain \beq \Theta^{E\ts
\tilde E}(u ,u)=\Theta^E_{i\bar j \alpha\bar\beta} u^{i\alpha
A}\bar{u^{j\beta B}} \cdot h^{\tilde E}_{A\bar B}+ \Theta^{\tilde
E}_{i\bar j A\bar B} u^{i\alpha A}\bar{u^{j\beta B}}\cdot
h^E_{\alpha\bar \beta}. \eeq It is nonnegative and one can see that
by  choosing  normal coordinates for $h^E$ and $h^{\tilde E}$ at a
fixed point. \eproof

\noindent We need  the following fundamental result in
\cite[Proposition~1.16]{DPS94}.

\blemma\label{key} Let $E\>X$ be a nef vector bundle over a compact
complex manifold $X$. Suppose $\sigma\in H^0(X,E^*)$ is a nonzero
section, then $\sigma$ does not vanish anywhere. \elemma

\noindent By refining the Bochner technique, we obtain the following
vanishing theorem for vector bundles with ``degenerate" curvature
tensors.

\btheorem \label{vanishing} Let $(X,\omega)$ be a compact K\"ahler
manifold of dimension $n$.  Let $(E,h)$ be a  holomorphic vector
bundle and $A$ be a line bundle.   Suppose either \bd\item $E$ is
nef and the second Ricci curvature $tr_\omega \Theta^E$ is
semi-positive, and   $A$ is  semi-ample but non-trivial; or
\item the second Ricci curvature $tr_\omega \Theta^E$ is strictly positive
and $A$ is nef. \ed
 Then \beq
H^{0}(X,E^*\ts A^*)=0.\eeq \etheorem

\bproof (1). Since $A$ is semi-ample, $A^k$ is generated by its
global sections for large $k$. Hence, there is an induced smooth
Hermitian metric $h^A$ on $A$ such that the curvature $\Theta^A$ is
semi-positive, i.e.
$$\sq\Theta^A=-\sq\p\bp\log h^A\geq 0.$$
On the other hand, since $A$ is not trivial, for the fixed K\"ahler
metric $\omega$ on $X$, the scalar curvature function $$tr_\omega
\Theta^A$$ is non-negative, but not identically zero. Indeed, if it
is identically zero, we deduce that $\Theta^A$ is identically zero,
and so $A$ is trivial.
 The curvature
tensor of $E\ts A$ can be written as $$\Theta^{E\ts A}=\Theta^E\ts
id_A+id_E\ts\Theta^A.$$ In order to prove (\ref{vanishing}), we
argue by contradiction. Suppose $H^0(X,E^*\ts A^*)\neq 0$, i.e.
there exists a nonzero section $\sigma\in H^0(X,E^*\ts A^*)$. By
Lemma \ref{key}, $\sigma$ is nowhere vanishing.  By using  standard
Bochner identity over $E^*\ts A^*$, \beq
\Delta_{\bar\partial}=\Delta_\p+[\sq\Theta^{E^*\ts
A^*},\Lambda_\omega],\eeq  we obtain \beq 0=\|\nabla
\sigma\|^2-(tr_\omega\Theta^{E^*\ts A^*}\sigma,\sigma). \eeq Note
that $\Theta^{E^*}=-(\Theta^E)^T$ as  $(1,1)$-form valued $r\times
r$ matrices.
 Hence,
\beq 0=\|\nabla
\sigma\|^2+\left(\left(\left[tr_\omega\Theta^{E}\right]^T\ts
id_A+id_E\ts tr_\omega\Theta ^A\right) \sigma,\sigma\right). \eeq By
using local holomorphic frames $\{z^i\}$,  $\{e_1,\cdots, e_r\}$,
$\{e\}$ of $X$, $E$ and $A$ respectively, we write $$\omega=\sq
g_{i\bar j}dz^i\wedge d\bar z^j, \ \Theta^E=R_{i\bar j
\alpha}^{\beta }dz^i\wedge d\bar z^j \ts e^\alpha \ts e_\beta,\ \ \
\Theta^A=R^A_{i\bar j}dz^i\wedge d\bar z^j,\ \ \sigma=\sigma^\alpha
e_\alpha\ts e.$$ We obtain \beq \int_X h^A\left(g^{i\bar j}R_{i\bar
j \beta\bar\alpha} \sigma^\alpha \bar \sigma^\beta +g^{i\bar
j}R^A_{i\bar j} h_{\alpha\bar \beta}\sigma^\alpha\bar\sigma^\beta
\right)\omega^n=0.\eeq  By assumption, the (transposed) second Ricci
curvature $$ \left[tr_\omega
\Theta^E\right]^T=\left(\sum_{i,j}g^{i\bar j} R_{i\bar j
\beta\bar\alpha}\right)$$ is Hermitian semi-positive as a $(r\times
r)$ matrix and so \beq g^{i\bar j}R_{i\bar j \beta\bar\alpha}
\sigma^\alpha \bar \sigma^\beta +g^{i\bar j}R^A_{i\bar j}
h_{\alpha\bar \beta}\sigma^\alpha\bar\sigma^\beta\equiv0.
\label{keyformula} \eeq It implies
$$g^{i\bar
j}R^A_{i\bar j} (h_{\alpha\bar
\beta}\sigma^\alpha\bar\sigma^\beta)\equiv0.$$ Since $\sigma$ is
nowhere vanishing,  we obtain $h_{\alpha\bar
\beta}\sigma^\alpha\bar\sigma^\beta>0$ at each point. Therefore
$$tr_\omega\Theta^A=g^{i\bar
j}R^A_{i\bar j}\equiv 0.$$ This is a contradiction.

(2). For  nef $A$, we use similar ideas as described  in the first
part of Lemma \ref{positive}.  Since $tr_\omega\Theta^E$ is strictly
positive,  there exists $\eps>0$ such that
$$(tr_\omega\Theta^E)^T(u(x),u(x))\geq (n+1)\eps |u(x)|_{ h}^2$$
for any $u\in\Gamma(X, E)$ since $X$ is compact. On the other hand,
 since $A$ is nef, there exists a smooth metric $h_0$ on the nef
line bundle $A$ such that $\sq \Theta^A \geq -\eps \omega_g. $
Hence, for any $u\in \Gamma(X, E)$ and $v\in \Gamma(X,A)$ \be
\left(\left[tr_\omega\Theta^{E}\right]^T\ts id_A+id_E\ts
tr_\omega\Theta ^A\right)(u\ts v,u\ts v)&\geq&
\left((tr_\omega\Theta^E)^T(u,u)-n\eps |u|_{ h}^2\right)
|v|^2_{h_0}\\&\geq& \eps |u|_{ h}^2|v|^2_{h_0}. \ee  Therefore in
(\ref{keyformula}), $\sigma\equiv0$, i.e. $H^0(X,E^*\ts A^*)=0$.
\eproof

\bremark The semi-positivity of the second Ricci curvature
$tr_\omega \Theta^E$ can be replaced by the semi-stability of $E$
with respect to $\omega$ following Donaldson-Uhlenbeck-Yau's
theorem. Moreover, the  Griffiths (or Nakano, or dual Nakano)
semi-positivity of $E$ can also imply the semi-positivity of the
second Ricci curvature $tr_\omega \Theta^E$.

 \eremark

 The following
Kawamata-Viehweg-Nadel type vanishing theorem for a semi-positive
vector bundle twisted by a big line bundle is essentially known to
experts (e.g., \cite[Theorem~4.2.4]{Ca98}, \cite[Theorem~5.11]{Dem},
\cite[Theorem~1.2]{Fuj12}, \cite[Theorem~1.1]{Fuj13},
\cite[Theorem~1.1]{Rau}, \cite{FM}, \cite{Ma}), although the
statement is not written down precisely. For the sake of
completeness, we include a short sketch here, following the approach
in \cite[Theorem~5.11]{Dem}  for line bundles.

\blemma\label{Nadel}  Let $(X,\omega)$ be a  K\"ahler weakly
pseudo-convex manifold, and let $A$ be a holomorphic line bundle
over $X$ equipped with a (possibly) singular Hermitian metric
$h=e^{-2\phi}$.  Assume that
$$\sq \Theta^{A}\geq \eps \omega$$ for some continuous positive
function $\eps$ on $X$. If $(E,h^E)$ is a Nakano semi-positive
vector bundle, then \beq H^q(X,K_X\ts E\ts A\ts \sI(\phi))=0 \eeq
for all $q\geq 1$. \elemma

\bproof Let $\mathscr{L}^q$ be the sheaf of germs of $(n,q)$-forms
$u$ with values in $E\ts A$ and with measurable coefficients such
that both $|u|_{h^E}^2\cdot e^{-2\phi}$ and $|\bp^{E\ts A} u|_{g\ts
h^E}\cdot  e^{-2\phi}$ are locally integrable.  The $\bp^{E\ts A}$
operator defines a  complex of sheaves $(\mathscr{L}^{\bullet},
\bp^{E\ts A})$ which is a fine resolution of the sheaf $\sO(K_X\ts
E\ts A)\ts \sI(\phi)$, i.e. we have the following exact sequence
\beq 0\>\sO(K_X\ts E\ts A)\ts \sI(\phi)\>\mathscr L^0\>\mathscr
L^1\>\cdots \>\mathscr L^n\>0.\eeq Indeed, it follows from a vector
bundle version of H\"ormander $L^2$-estimate
(\cite[Corollary~5.3]{Dem}) since the vector bundle $E\ts A$ has a
singular metric which is Nakano positive in the sense of current. By
using the $L^2$ estimate again (e.g. \cite[Theorem~5.11]{Dem}), one
can show $H^q(\Gamma(X,\mathscr L^\bullet))=0$ for $q\geq 1$ and we
obtain the desired vanishing cohomologies. \eproof

\noindent Next, we introduce two different concepts on numerical
dimension for nef and pseudo-effective line bundles.

 \bdefinition Let $N$ be a \emph{nef
line bundle} over a compact K\"ahler manifold $X$ with $\dim_\C
X=n$. The numerical dimension $\nu(N)$ of $N$ is defined as
(\cite[Definition~6.20]{Dem}) \beq \nu(N)=\max\{k=0,\cdots,n\ | \
c^k_1(N)\neq 0 \in \text{$H^{2k}(X,\R)$}.\}\eeq \edefinition

\bdefinition

 Let $(A,e^{-2\phi})$ be a \emph{pseudo-effective line bundle} over a
 compact K\"ahler manifold $X$ of dimension $n$. The \emph{numerical dimension} $\text{nd}(A,\phi)$
 of $A$ is defined   in \cite[Definition~3.1]{Cao14} as the largest number such that the cohomological product $\la (\sq \p\bp\phi)^k\ra\neq 0$.

\edefinition

  Note that, in general these two definitions do not coincide even for nef line bundles (\cite[Remark~7]{Cao14}, or \cite[Remark~4.3.5]{Dem14}).
 For a discussion of the relationship between various definitions of numerical
 dimensions, we refer to the paper \cite[Section~4.3]{Dem14}.\\

 Based on their solution to Demailly's strong openness conjecture
(\cite{GZ15}), Guan-Zhou achieved in \cite[Corollary~1.7]{GZ} (see
also \cite[Corollary~3.2]{GZ15}) a celebrated Kawamata-Viehweg-Nadel
type vanishing theorem which generalizes a theorem in
\cite[Theorem~1.3]{Cao14} (see also Demailly's survey paper
\cite{Dem14} for a variety of vanishing theorems).

 \blemma\label{GZ} Let $(A,e^{-2\phi})$ be a {pseudo-effective line bundle} over a
 compact K\"ahler manifold $X$ of dimension $n$. Then for any $q>n-\emph{nd}(A,\phi)$,
 \beq H^q(X,K_X\ts A\ts \sI(\phi))=0.\eeq
 \elemma

\vskip 1\baselineskip

\noindent We shall use the ideas in the proof of Lemma \ref{GZ}
(e.g. \cite{Cao14} and \cite{GZ15}) and Lemma \ref{Nadel} to prove
Theorem \ref{mainkey}.

\emph{The proof of Theorem \ref{mainkey}.} Suppose
$\text{nd}(A,\phi)=n$. By using Guan-Zhou's solution to Demailly's
strong openness conjecture (\cite{GZ15}) and the construction in
\cite[Lemma~5.5]{Cao14}, there exists a singular metric
$e^{-2\phi_1}$ on the pseudo-effective line bundle $A$ such that
$$\sI(\phi_1)=\sI(\phi)$$
and it is curvature current $$\sq\p\bp\phi_1>c\omega$$ for some
smooth positive $(1,1)$ form $\omega$ and constant $c>0$. By
applying Lemma \ref{Nadel} to $(X,E,A, \sI(\phi_1))$, we get \beq
H^q(X,K_X\ts E\ts A\ts \sI(\phi))=0, \qtq{for} q>0.\eeq

 Now we assume that $\text{nd}(A,\phi)<n$. We use similar ideas as in \cite[Theorem~6.25]{Dem}, \cite[Proposition~5.6]{Cao14}.  Let $B$
be a very  ample divisor such that $B\ts A$ is ample and
$\sI(\phi|_B)=\sI(\phi)|_B$ (\cite[Theorem~1.10]{FM}). We consider
the exact sequence
  $$0\>\sO_X(-B)\>\sO_X\>(i_B)_*\sO_B\>0.$$
By tensoring with $K_X\ts E\ts A\ts B\ts  \sI(\phi)$ and using
adjunction formula, one gets
\be &&\cdots\>H^q(X,K_X\ts E\ts A\ts \sI(\phi))\>H^q(X,K_X\ts E\ts A\ts B\ts \sI(\phi))\\
&&\>H^q(B,K_B\ts (E\ts A)|_B\ts\sI(\phi|_B) )\>H^{q+1}(X, K_X\ts
E\ts A\ts\sI(\phi))\>\cdots.\ee Since $E$ is Nakano semi-positive,
and $A\ts B$ is ample, by Lemma \ref{Nadel}, we have \beq
H^q(X,K_X\ts E\ts A\ts B\ts \sI(\phi))=0\qtq{for }q>0. \eeq
Therefore \beq H^q(B,K_B\ts (E\ts A)|_B\ts\sI(\phi|_B) )\cong
H^{q+1}(X, K_X\ts E\ts A\ts\sI(\phi)) \label{key111}\eeq for every
$0<q<n$. Moreover, {$E|_B$ is also Nakano semi-positive} over $B$.
Hence the induction hypothesis implies that the cohomology group on
$B$ on the right hand side of (\ref{key111}) is zero when
$q>n-\text{nd}(A,\phi)$. \qed

\vskip 1\baselineskip

Similarly, we  have the following  variant of Theorem \ref{mainkey}
(see also \cite[Theorem~6.17]{Dem} for the line bundle case.)
\bproposition\label{nakanoKV} Let $X$ be a smooth projective
variety. Let $A$ be a line bundle over $X$ such that some positive
multiple $mA$ can be written as $mA=N+D$ where $N$ is a nef line
bundle and $D$ is an effective divisor. If $(E,h)$ is a Nakano
semi-positive vector bundle, then \beq H^q(X,K_X\ts E\ts A\ts
\sI(m^{-1}D))=0 \qtq{for} q>n-\nu(N), \eeq where $\nu(N)$ is the
numerical dimension of the nef line bundle $N$.
 \eproposition

 As a special case, one has

\bcorollary \label{nakanoKV1} If $(E,h^E)$ is a Nanako semi-positive
vector bundle and $A$ is a nef line bundle, then \beq H^q(X,E^*\ts
A^*)=0 \qtq{if} q<\nu(A).\label{s1}\eeq In particular, if in
addition, $A$ is not numerically trivial, then \beq H^0(X,E^*\ts
A^*)=0.\eeq \ecorollary

\bproof We can take $m=1$, $D=0$ and $A=N$ in Proposition
\ref{nakanoKV}. By Serre duality, we obtain (\ref{s1}). Hence, when
$\nu(A)\geq 1$, or equivalently, $\nu(A)\neq 0$, $H^0(X,E^*\ts
A^*)=0.$ It is well-known that for a nef line bundle $A$, $\nu(A)=0$
if and only if $A$ is numerically trivial. \eproof

\vskip 1\baselineskip

 As an application of Theorem \ref{mainkey}, we obtain
general vanishing theorems on the relative setting which generalize
Koll{\'a}r's vanishing theorems.

\btheorem\label{relativekey} Let $f:X\>Y$ be a holomorphic
submersion between two smooth complex projective varieties and
$\dim_\C Y=n$. Let $(A,e^{-2\phi})$ be a pseudo-effective line
bundle over $X$.  If both $E_1$ and $E_2$ are  Nakano semi-positive
vector bundles, then for any $s\geq 1$ and
$q>n-\emph{nd}(A,\phi)$, we have
$$ H^q\left(Y,f_*(K_{X/Y}\ts E_1)^{\ts s}\ts K_Y\ts E_2\ts A\ts
\sI(\phi)\right)=0, $$ as long as $f_*(K_{X/Y}\ts E_1)$ is locally
free.
 \etheorem

\bproof By using Theorem \ref{mainkey}, we only need to show that
$$f_*(K_{X/Y}\ts E_1)^{\ts s}\ts E_2$$ is Nakano semi-positive for any
$s\geq 1$. Indeed by a result of \cite{MT08} (see also \cite{LY15}),
if  $f_*(K_{X/Y}\ts E_1)$ is locally free, then $f_*(K_{X/Y}\ts
E_1)$ has a Nakano semi-positive metric.  By part $(3)$ of Lemma
\ref{positive}, we deduce $$f_*(K_{X/Y}\ts E_1)^{\ts s}\ts E_2$$ is
Nakano semi-positive for any $s\geq 1$. \eproof

\bremark Theorem \ref{relativekey} also holds when $X$ and $Y$ are
compact K\"ahler manifolds. \eremark

\vskip 2\baselineskip

\section{Fujita Conjecture type theorems on complex projective varieties}

In this section, we derive Fujita Conjecture type theorems on
complex projective varieties and prove Theorem \ref{main11}, Theorem
\ref{main2},  Theorem \ref{main3} and Theorem \ref{main4}.\\

 \emph{The proof of Theorem \ref{main11}.}  We show the coherent
 sheaf $K_X\ts L^{\ts n}\ts A\ts \sI(\phi)$ is $0$-regular, and so by Lemma \ref{cmr}, it is globally generated. When
 $0<q<n$, $L^{\ts (n-q)}\ts A$ is indeed a big line bundle, and we
 can apply
 Nadel vanishing theorem (\cite{Nad90} or \cite[Theorem~0.3]{Dem14} or
 Lemma
 \ref{GZ}),
 \beq H^q(X,K_X\ts L^{\ts (n-q)}\ts A\ts \sI(\phi))=0.\eeq
When $q=n$, we need to show  \beq H^n(X,K_X\ts A\ts
\sI(\phi))=0,\eeq which follows from  Lemma \ref{GZ} and the
assumption  $\text{nd}(A,\phi)\geq 1$. \qed

\vskip 1\baselineskip

\emph{The proof of Theorem \ref{main2}.}  By Castelnuovo-Mumford
regularity (e.g. Lemma \ref{cmr}), we only need to prove $K_X \ts
L^{\ts n}\ts \left(E\ts A\right)$ is $0$-regular with respect to
$L$. Hence, it suffices to show \beq H^q(X,K_X \ts L^{\ts (n-q)}\ts
(E\ts A))=0 \qtq{for all} q>0.\eeq

For $0<q<n$, we claim that the vector bundle $L^{\ts (n-q)}\ts (E\ts
A)$ has a smooth metric with strictly positive curvature in the
sense of Nakano, and by Nakano vanishing theorem( Lemma
\ref{nakano}), we have the desired vanishing cohomologies. Indeed,
since $n-q\geq 1$, $L^{\ts(n-q)}\ts A$ is ample. By  Lemma
\ref{positive}, $E\ts L^{\ts (n-q)}\ts A$ is strictly positive in
the sense of Nakano.

 When $q=n$, we need to show
$H^n(X,K_X \ts E\ts A)=0$ or equivalently $H^0(X,E^*\ts A^*)=0$ when
$E$ is Nakano semi-positive and $A$ is nef but not numerically
trivial. This is assured by Corollary \ref{nakanoKV1}.

 The case when $q>n$ is obvious and  we complete the proof of
 Theorem \ref{main2}. \qed

\vskip 1\baselineskip

 We have the following variant of Theorem \ref{main2}.
\btheorem\label{main21} Let $X$ be a compact K\"ahler manifold and
$L$ be an ample and globally generated line bundle. Let $(E,h)$ be a
Hermitian  holomorphic vector bundle with Nakano positive curvature.
Suppose $A$ to be a nef line bundle, then the vector bundle \beq K_X
\ts L^{\ts n} \otimes E\ts A\eeq is globally generated. \etheorem

\bproof We use similar ideas as described in the proof of Theorem
\ref{main2}. Suppose $E$ is Nakano positive and $A$ is nef,  then by
Lemma \ref{positive}, the vector bundle $E\ts A$ admits a smooth
Hermitian metric whose curvature is strictly positive in the sense
of Nakano. In particular, the second Ricci curvature
$tr_\omega\Theta^{E\ts A}$ is strictly positive. We obtain the
vanishing cohomology in Theorem \ref{vanishing}. Finally, one can
follow the steps in the proof of Theorem \ref{main2}. \eproof

\noindent It is proved in  \cite[Theorem~1.2]{Ber09} that if $E$ is
a semi-ample (resp. ample) vector bundle, then $E\ts \det E$ is
Nakano semi-positive (resp. Nakano positive). Hence, by Theorem
\ref{main2} and Theorem \ref{main21}, we get

\bcorollary\label{algebraic} Let $X$ be a smooth complex projective
variety and $L$ be an ample and globally generated line bundle. Let
$E$ be a vector bundle and $A$ be a line bundle. Suppose either \bd
\item $E$ is  semi-ample, $A$ is nef but not numerically trivial; or

\item $E$ is ample and $A$ is nef.

\ed Then the vector bundle \beq K_X \ts L^{\ts n} \otimes (E\ts \det
E)\ts A\eeq is globally generated. \ecorollary

Thanks to Theorem \ref{mainkey}, one can also get the following
variant of Theorem \ref{main2}. \btheorem\label{main211} Let
$(X,\omega)$ be a compact K\"ahler manifold and $L\>X$ be an ample
and globally generated line bundle. Suppose $(A,e^{-2\phi})$ is a
pseudo-effective line bundle and $\sI(\phi)$ is the multiplier ideal
sheaf. Let $E$ be a Nakano semi-positive vector bundle.  If the
numerical dimension $\emph{nd}(A,\phi)\neq 0$, then
 $$K_X \ts L^{\ts n}\ts E \otimes A\ts \sI(\phi)$$ is
globally generated. \etheorem

\vskip 1\baselineskip

\emph{The proof of Theorem \ref{main3}.} By the assumption, the
direct image sheaf of the relative canonical line bundle
$f_*(K_{X/Y})$ is indeed a holomorphic vector bundle.   In the
literatures, it is known that the vector bundle  $f_*(K_{X/Y})$ is
weakly positive in suitable  sense (e.g. \cite[Theorem~III]{Vie83},
\cite[Corollary~5]{GT84}, \cite[Corollary~3.7]{Kol86}). Here, we use
a recent fact \cite[Theorem~1.2]{Ber09} of Berndtsson that
$f_*(K_{X/Y})$ is actually semi-positive in the sense of Nakano. Let
$E=f_*(K_{X/Y})^{\ts s}$ and so \beq f_*(K_{X/Y})^{\ts s}\ts K_Y\ts
L^{\ts n}\ts A= K_Y\ts L^{\ts n}\ts (E\ts A).\eeq  According to
Lemma \ref{positive}, $E$ is also semi-positive in the sense of
Nakano. By Theorem \ref{main2}, $K_Y\ts L^{\ts n}\ts (E\ts A)$ is
globally generated as long as $A$ is nef but not numerically
trivial. The proof of Theorem \ref{main3} is completed. \qed

\vskip 1\baselineskip

As an application of Theorem \ref{relativekey}, we have the
following slightly general version of Theorem \ref{main3} .
\btheorem\label{main311} Let $f:X\>Y$ be a holomorphic submersion
between two smooth complex projective varieties and $\dim_\C Y=n$.
Suppose that $L\>Y$ is an ample and globally generated line bundle, and
$(A,e^{-2\phi})\>Y$ be a pseudo-effective line bundle with
$\emph{nd}(A,\phi)\neq 0$. Suppose both $E_1$ and $E_2$ are Nakano
semi-positive, then for any $s\geq 1$ \beq f_*(K_{X/Y}\ts E_1)^{\ts
s}\ts K_Y \ts L^{\ts n} \otimes E_2\ts A\ts \sI(\phi)\eeq is
globally generated as long as $f_*(K_{X/Y}\ts E_1)$ is locally free.
\etheorem

\vskip 1\baselineskip \emph{The proof of Theorem \ref{main4}.} When
the family $\sX\>S$ is effectively parameterized, Schumacher proved
in \cite[Theorem~1]{Schu} that the naturally induced Hermitian
metric on the relative canonical line bundle $K_{\sX/S}$ is
 strictly positive. By
\cite[Corollary~2]{Schu} or \cite[Theorem~1.2]{Ber09}, we know \beq
f_*(K^{\ts s}_{X/S})\eeq is strictly Nakano positive for all $s>1$.
One can write \beq f_*(K^{\ts s}_X)\ts L^{\ts n}\ts A =f_*(K^{\ts
s}_{X/S})\ts K_S \ts L^{\ts n}\ts (K^{\ts (s-1)}_S\ts A)\eeq We
first observe that the canonical bundle $K_S$ is  nef. Indeed, by a
recent result in \cite[Theorem~1]{TY15}, the compact complex base
$S$ is actually (Kobayashi) hyperbolic. Hence, it contains no
rational curve. By using Mori's cone theorem \cite{Mori82}, we
deduce $K_S$ is nef since it is a projective manifold without
rational curve. The global generation of
 the vector bundle $f_*(K^{\ts s}_{X/S})\ts K_S
\ts L^{\ts n}\ts (K^{\ts (s-1)}_S\ts A)$  follows from Theorem
\ref{main21} since $f_*(K^{\ts s}_{X/S})$ is Nakano positive and
$K^{\ts (s-1)}_S\ts A$ is nef. \qed

\bremark On smooth complex projective varieties, we can also derive
globally generation for symmetric powers and wedge powers of vector
bundles by using the corresponding vanishing theorems (e.g.
\cite{LY15}). \eremark

\begin{center}Acknowledgements.\end{center}  This work was partially supported by
China's Recruitment
 Program of Global Experts and NSFC 11688101. It was carried out when the
authors attended  the seminar on ``Fujita conjecture and related
topics
  " held at Math Institute of Chinese Academy of Sciences, and the authors would like to thank all members of this seminar for  stimulating discussions.
 In particular, the authors  wish to thank  Yifei Chen,  Baohua Fu and  Xiaotao
 Sun  for  some useful suggestions. The second
 author  would like to thank Professor S.-T. Yau
 for his advice, support and encouragement. He also
 thanks  Kefeng Liu, Valentino Tosatti and Xiangyu Zhou
 for  helpful discussions. The authors would also like to thank the anonymous
referees whose comments and suggestions helped improve and clarify
the paper.

\vskip 1\baselineskip

\end{document}